\documentclass[leqno, 12pt]{article}
\usepackage{amsmath,amsfonts,amsthm,amssymb,indentfirst}

\newtheorem{theorem}{Theorem}
\newtheorem{lemma}[theorem]{Lemma}
\newtheorem{corollary}[theorem]{Corollary}
\newtheorem{proposition}[theorem]{Proposition}

\newtheorem{definition}[theorem]{Definition}

\newcommand{\E}{\mathrm{End}}
\newcommand{\preq}{\preccurlyeq}

\begin{document}

\title{Endomorphism rings generated using small numbers of elements}
\author{Zachary Mesyan}

\maketitle

\begin{abstract}
Let $R$ be a ring, $M$ a nonzero left $R$-module, and $\Omega$ an
infinite set, and set $E = \E_R (\bigoplus_{\Omega} M)$. Given two
subrings $S_1, S_2 \subseteq E$, write $S_1 \approx S_2$ if there
exists a finite subset $U \subseteq E$ such that $\langle S_1 \cup
U \rangle = \langle S_2 \cup U \rangle$. We show that if $M$ is
simple and $\Omega$ is countable, then the subrings of $E$ that
are closed in the function topology and contain the diagonal
subring of $E$ (consisting of endomorphisms that take each copy of
$M$ to itself) fall into exactly two equivalence classes, with
respect to the equivalence relation above. We also show that every
countable subset of $E$ is contained in a $2$-generator
subsemigroup of $E$.
\end{abstract}

\section{Introduction}

Let $R$ be a ring, $M$ a nonzero left $R$-module, $\Omega$ an
infinite set, $N = \bigoplus_\Omega M$, and $E = \E_R (N)$. In
this paper we will show that the ring $E$ has some unusual
properties that are analogous to known properties of the symmetric
group of an infinite set.  In the next section we will demonstrate
that every countable subset of $E$ is contained in a 2-generator
subring of $E$. (The proof also works if $N$ is taken to be
$\prod_\Omega M$.) This parallels Galvin's result that every
countable subset of the symmetric group of an infinite set is
contained in a $2$-generator subgroup (cf.~\cite{Galvin}). As an
immediate corollary, we will show that every countable ring can be
embedded in a ring generated by two elements, reproducing a result
of Maltsev (cf.~\cite{O'Meara} and~\cite{OVW} for different
proofs). The group-theoretic analog of this fact has also been
known for a long time (cf.~\cite{Higman} and~\cite{Neumann}).

Actually, our proof of the above result shows that a countable
subset of $E$ is contained in a $2$-generator {\em subsemigroup}
of $E$. This is a generalization of the result of Magill
in~\cite{Magill} that every countable set of endomorphisms of an
infinite-dimensional vector space is contained in a 2-generator
subsemigroup of the semigroup of all endomorphisms of that vector
space (see also~\cite{AMS} for a shorter proof).

Given two subrings $S_1, S_2 \subseteq E$, we will say that $S_1
\approx S_2$ if there exists a finite subset $U \subseteq E$ such
that $\langle S_1 \cup U \rangle = \langle S_2 \cup U \rangle$. We
will devote the remainder of the paper to exploring properties of
this equivalence relation. In particular, we will show that if $M$
is finitely generated, $\E_R(M)$ is a simple ring (e.g., if $M$ is
a simple module), $\Omega$ is countable, and $S \subseteq E$ is a
subring that is closed in the function topology and contains $D$,
the diagonal subring of $E$ (consisting of endomorphisms that take
each copy of $M$ to itself), then either $S \approx D$ or $S
\approx E$ (but not both). This is in the spirit the result of
Bergman and Shelah in~\cite{B&S} that the subgroups of the
symmetric group of a countably infinite set that are closed in the
function topology fall into exactly four equivalence classes.
(There the equivalence relation is defined the same way as the
relation above, with subgroups in place of subrings.)

Along the way, we will also note a natural way of associating to
every preordering $\rho$ on $\Omega$ a subring $E(\rho)$ of $E$.
We will then show that if $\Omega$ is countable and $M$ is
finitely generated, then the subrings of the form $E(\rho)$ fall
into exactly two equivalence classes (again, represented by $D$
and $E$). The result mentioned in the previous paragraph is
actually a special case of this.

A curious example is that if we view $E$ as a ring of row-finite
matrices over $\E_R(M)$, then the subring of upper-triangular
matrices is equivalent to $E$, while the subring of
lower-triangular matrices is equivalent to $D$.

I am grateful to George Bergman, whose comments and suggestions
have led to vast improvements in this paper, and also to Kenneth
Goodearl for referring me to related literature.

\section{Countable sets of endomorphisms}

Let $R$ denote a unital associative ring, $M$ a nonzero left
$R$-module, and $\Omega$ an infinite set. $N$ will denote either
$\bigoplus_{\alpha \in \Omega} M_\alpha$ or $\prod_{\alpha \in
\Omega} M_\alpha$ (the arguments in this section work under either
interpretation), where each $M_\alpha = M$. We will write $E$ to
denote $\E_R (N)$. Endomorphisms will be written on the right of
their arguments. Also, given a subset $\Sigma \subseteq \Omega$,
we will write $M^{\Sigma}$ for the $R$-submodule of $N$ consisting
of elements $(n_\alpha)_{\alpha \in \Omega}$ with $n_\alpha = 0$
for all $\alpha \notin \Sigma$, and $\pi_{\Sigma}$ for the
projection from $N$ to $M^{\Sigma}$ along
$M^{\Omega\backslash\Sigma}$, so in particular, $N = M^{\Sigma}
\oplus M^{\Omega\backslash\Sigma}$. Finally, $\mathbb{Z}^+$ will
denote the set of positive integers, and if $\Gamma$ is a set,
$|\Gamma|$ will denote the cardinality of $\Gamma$.

The following argument was obtained by tinkering with the proofs
of Theorem 2.6 in~\cite{Magill} and Theorem 3.1 in~\cite{Galvin}.

\begin{theorem}\label{4-gen}
Every countable subset of $E = \E_R (N)$ is contained in a
$2$-generator subsemigroup of $E$ $($viewed as a multiplicative
semigroup$)$.
\end{theorem}

\begin{proof}
We may assume that $\Omega = \mathbb{Z} \times \Gamma$, where
$|\Gamma| = |\Omega|$.  Let us set $\Sigma = \{0\} \times \Gamma$.
Also, let $g_1 \in E$ be an endomorphism that takes $N$
isomorphically to $M^\Sigma$, and let $g_2 \in E$ be the right
inverse of $g_1$ that takes $M^\Sigma$ isomorphically to $N$ and
takes $M^{\Omega \setminus \Sigma}$ to zero. We note that $g_2g_1
= \pi_\Sigma$.

Now, let $U \subseteq E$ be a countably infinite subset.  We will
show that $U$ is contained in a 2-generator subsemigroup of $E$.
Since $E = g_1 \pi_{\Sigma} E \pi_{\Sigma} g_2$, we can find a
subset $\bar{U} \subseteq \pi_{\Sigma} E \pi_{\Sigma}$ such that
$U = g_1 \bar{U} g_2$. Since $\bar{U}$ is countable, we can write
$\bar{U} = \{u_i : i \in \mathbb{Z}\}$. For each $i \in
\mathbb{Z}$ let us define $\hat{u_i} \in \E_R (M^{\{i\} \times
\Gamma})$ so that $\hat{u_i}$ acts on $M^{\{i\} \times \Gamma}$ as
$u_i$ acts on $M^{\Sigma}$ (upon identifying $M_{(i,\ \gamma)}$
with $M_{(0,\ \gamma)}$ for each $\gamma \in \Gamma$). Also, let
$g_3 \in E$ be an endomorphism such that for each $i \in
\mathbb{Z}$ the restriction of $g_3$ to $M^{\{i\} \times \Gamma}$
is $\hat{u_i}$. Finally, let $g_4 \in E$ be the automorphism that
takes $M_{(i,\ \gamma)}$ identically to $M_{(i+1,\ \gamma)}$ for
each $i \in \mathbb{Z}$ and $\gamma \in \Gamma$. Then for each $i
\in \mathbb{Z}$, we have $u_i = \pi_{\Sigma} g_4^i g_3 g_4^{-i} =
g_2g_1g_4^i g_3 g_4^{-i}$. Writing $g_5 = g_4^{-1}$ and recalling
that $g_2$ is a right inverse of $g_1$, we conclude that $U = \{
g_1 g_4^i g_3 g_5^i g_2 : i \in \mathbb{Z} \}$.

It remains to be shown that $\{g_1, g_2, g_3, g_4, g_5\}$ is
contained in a $2$-generator subsemigroup of $E$. Write $\Omega$
as $\bigcup_{i=1}^7 \Sigma_i$, where the union is disjoint, and
$|\Sigma_i| = |\Omega|$ for $i \in \{1, 2, \dots, 7\}$. Also,
write $\Delta = \bigcup_{i=2}^7 \Sigma_i$. Now, let us choose an
endomorphism $f_1 \in E$ such that
\begin{enumerate}
\item[$(1)$] $f_1$ takes $M^{\Sigma_i}$ isomorphically to
$M^{\Sigma_{i+1}}$ for $i \in \{1, 2, \dots, 5\}$, and takes
$M^{\Sigma_6 \cup \Sigma_7}$ isomorphically to $M^{\Sigma_7}$.
\end{enumerate}
Such an endomorphism will necessarily map $N$ isomorphically to
$M^\Delta$. Let $f_2 \in E$ be an endomorphism such that
\begin{enumerate}
\item[$(2)$] $f_2$ takes $N$ isomorphically to $M^{\Sigma_1}$.
\end{enumerate}
Also, let us define $t_i \in \mathrm{Hom}_R (M^{\Sigma_1},
M^{\Sigma_{i+1}})$ ($i \in \{1, 2, \dots, 5\}$) so that
\begin{enumerate}
\item[$(3)$] $t_i$ is the restriction of $f_1^i$ to $M^{\Sigma_1}$.
\end{enumerate}
Then $t_i^{-1}f_2^{-1}g_i \in \mathrm{Hom}_R (M^{\Sigma_{i+1}},
N)$ for $i \in \{1, 2, \dots, 5\}$, where $t_i^{-1}, f_2^{-1} \in
E$ are right inverses of $t_i$ and $f_2$, respectively. Writing
\begin{enumerate}
\item[$(4)$] $t_6 = f_1^6$,
\end{enumerate}
we see that $t_6$ is an isomorphism from $N$ to $M^{\Sigma_7}$ and
$t_6^{-1}f_2 \in \mathrm{Hom}_R (M^{\Sigma_7}, N)$. Finally, let
$f_3 \in E$ be an endomorphism such that
\begin{enumerate}
\item[$(5)$] $f_3$ restricted to $M^{\Sigma_{i+1}}$ is
$t_i^{-1}f_2^{-1}g_i$ if $i \in \{1, 2, \dots, 5\}$, and
$t_6^{-1}f_2$ if $i = 6$.
\end{enumerate}
Then $g_i = f_1^6f_3f_1^if_3$ for $i \in \{1, 2, \dots, 5\}$, and
therefore $\{g_1, g_2, g_3, g_4, g_5\}$ is contained in the
subsemigroup generated by $f_1$ and $f_3$.
\end{proof}

The following was originally proved by Maltsev.

\begin{corollary}
Every countable ring can be embedded in a ring generated by two
elements, using an embedding that respects central elements.
\end{corollary}

\begin{proof}
Let $S$ be a countable ring and $\Omega$ an infinite set. Then $S$
embeds diagonally in $\E_S (\bigoplus_\Omega S)$. Thus, by the
previous theorem, the image of $S$ is contained in a $2$-generator
subring of $\E_S (\bigoplus_\Omega S)$. It is clear that the
diagonal embedding maps the center of $S$ into the center of $\E_S
(\bigoplus_\Omega S)$.
\end{proof}

\section{Equivalence classes}

In this section we will keep $R$, $M$, $\Omega$, and $E$ as above,
but restrict our attention to the case $N = \bigoplus_\Omega M$.
However, we begin with two definitions applicable to an arbitrary
ring.

\begin{definition}
Let $S$ be a ring, $\kappa$ an infinite cardinal, and $S_1,\ S_2$
subrings of $S$.  We will write $S_1 \preq_{\kappa,S} S_2$ if
there exists a subset $U \subseteq S$ of cardinality $<\kappa$
such that $S_1 \subseteq \langle S_2 \cup U \rangle$, the subring
of $S$ generated by $S_2 \cup U$. If $S_1 \preq_{\kappa,S} S_2$
and $S_2 \preq_{\kappa,S} S_1,$ we will write $S_1
\approx_{\kappa,S} S_2$, while if $S_1 \preq_{\kappa,S} S_2$ and
$S_2 \not\preq_{\kappa,S} S_1$, we will write $S_1
\prec_{\kappa,S} S_2$.  The subscripts $S$ and $\kappa$ will be
omitted when their values are clear from the context.
\end{definition}

It is easy to see that $\preq_{\kappa,S}$ is a preorder on
subrings of $S$, and hence $\approx_{\kappa,S}$ is an equivalence
relation. This equivalence relation and the results of this
section are modeled on those in~\cite{B&S}, where Bergman and
Shelah define an analogous relation for groups and classify into
equivalence classes the subgroups of the group of permutations of
a countably infinite set that are closed in the function topology.
Properties of such an equivalence relation defined for submonoids
of the monoid of self-maps of an infinite set are investigated
in~\cite{ZM2}.

\begin{definition}
Let $S$ be a ring. Then the {\em cofinality} $c(S)$ of $S$ is the
least cardinal $\kappa$ such that $S$ can be expressed as the
union of an increasing chain of $\kappa$ proper subrings.
\end{definition}

Cofinality can be defined analogously for any algebra (in the
sense of universal algebra). It has received much attention in the
literature in connection with permutation groups. In particular,
Macpherson and Neumann show in~\cite{DM&PN} that
$c(\mathrm{Sym}(\Omega)) > |\Omega|$, where $\mathrm{Sym}(\Omega)$
is the group of all permutations of an infinite set $\Omega$. It
is shown in~\cite{ZM} that the ring $E$ likewise satisfies $c(E) >
|\Omega|$.

\begin{proposition}
Let $S, S' \subseteq E$ be subrings.\\[2pt]
{\rm(i)} $S \preq_{\aleph_0} S'$ if and only if $S
\preq_{\aleph_1} S'$ {\rm(}and hence $S \approx_{\aleph_0} S'$
if and only if $S \approx_{\aleph_1} S')$.\\[2pt]
{\rm(ii)} $S \approx_{\aleph_0} E$ if and only if $S
\approx_{|\Omega|^+} E$ $($where $|\Omega|^+$ is the successor
cardinal of $|\Omega|$$)$.
\end{proposition}

\begin{proof}
(i) follows from Theorem~\ref{4-gen}.  (ii) follows from the fact
that $c(E) > |\Omega|$. For, if $S \approx_{|\Omega|^+} E$, then
among subsets $U \subseteq E$ of cardinality $\leq |\Omega|$ such
that $\langle S \cup U \rangle = E$, we can choose one of least
cardinality. Let us write $U = \{f_i : i \in |U|\}$. Then the
subrings $S_i = \langle S \cup \{f_j : j < i \} \rangle$ $(i \in
|U|)$ form a chain of $\leq |\Omega|$ proper subrings of $E$.  If
$|U|$ were infinite, this chain would have union $E$,
contradicting $c(E) > |\Omega|$. Hence, $U$ is finite, and $S
\approx_{\aleph_0} E$.
\end{proof}

We will devote the rest of this section to showing that a large
natural class of subrings of $E$ consists of elements that are
$\prec_\kappa E$. First, we need a few definitions and a lemma.

\begin{definition}
Let $S$ be a ring and $U$ a subset of $S$.  We will say that $s
\in S$ is represented by a \emph{ring word of length 1} in
elements of $U$ if $r \in U \cup \{0, 1, -1\}$, and, recursively,
that $s \in S$ is represented by a \emph{ring word of length} $n$
in elements of $U$ if $s = p + q$ or $s = pq$ for some elements
$p, q \in S$ which can be represented by ring words of lengths
$m_1$ and $m_2$ respectively, with  $n = m_1 + m_2$.
\end{definition}

\begin{definition}
Let $U \subseteq E$ be a subset and $x_1, x_2 \in N =
\bigoplus_{\alpha \in \Omega} M_\alpha$. We will write $p_U (x_1,
x_2) = r$ if $x_2 = x_1 f$ for some $f \in E$ that is represented
by a ring word of length $r$ in elements of $U$, and $r$ is the
smallest such integer.  If no such integer exists, we will write
$p_U (x_1, x_2) = \infty$.  Also, given $x \in N$ and $r \in
\mathbb{Z}^+$, let $B_U (x, r) = \{y \in N : p_U (x, y) \leq r
\}$.  $\mathrm{(}$Here $p$ stands for ``proximity," and $B$ stands
for ``ball."$\mathrm{)}$
\end{definition}

\begin{definition}
Given a nonzero subset $X \subseteq N$, we will say that $\Sigma
\subseteq \Omega$ is the \emph{support} of $X$ if $X \subseteq
M^\Sigma$ and $\Sigma$ is the least such subset of $\Omega$. Also,
if $\kappa$ is a regular infinite cardinal $\leq |\Omega|$, we
will say that a subring $S \subseteq E$ is $\kappa$-\emph{fearing}
if for every $\alpha \in \Omega$, $(M_\alpha) S$ has support of
cardinality $< \kappa$.
\end{definition}

\begin{lemma}\label{bowls}
Suppose that $\kappa$ is a regular infinite cardinal $\leq
|\Omega|$, that $M$ can be generated by $< \kappa$ elements as an
$R$-module, that $S \subseteq E$ is a $\kappa$-fearing subring,
and that $U \subseteq E$ is a subset of cardinality $< \kappa$.
Then for any $x \in N$ and any $r \in \mathbb{Z}^+$, $B_{S \cup U}
(x, r)$ has support of cardinality $< \kappa$.
\end{lemma}

\begin{proof}
Let $x \in N$ be any element.  Then $|\{ xf : f \in U \}| <
\kappa$, since $|U| < \kappa$.  Hence, the support of $\{ xf : f
\in U \}$ is contained the union of $< \kappa$ finite sets and
therefore has cardinality $< \kappa$. Also, $\{ xf : f \in S \}$
has support of cardinality $< \kappa$. (There is a finite subset
$\Gamma \subseteq \Omega$ such that $x \in M^\Gamma$. So, since
$S$ is $\kappa$-fearing, $(M^\Gamma) S$ has support of cardinality
$< \kappa$.)  Therefore, $B_{S \cup U} (x, 1) = \{ xf : f \in S
\cup U \}$ has support of cardinality $< \kappa$.

Now, let $X \subseteq N$ be a subset that has support of
cardinality $< \kappa$.  Since $M$ can be generated by $< \kappa$
elements, $X$ is contained in a submodule of $N$ that can be
generated by $< \kappa$ elements.  Let $\{x_\varphi : \varphi \in
\Phi \}$ be a generating set for such a submodule, where $|\Phi| <
\kappa$. Then the submodule of $N$ generated by $\{ Xf : f \in S
\cup U \}$ is contained in the submodule of $N$ generated by
$\bigcup_{\varphi \in \Phi} \{x_\varphi f : f \in S \cup U \}$,
which has support of cardinality $< \kappa$, by the previous
paragraph and the regularity of $\kappa$.  Thus, $\{ Xf : f \in S
\cup U \}$ has support of cardinality $< \kappa$ as well.

Hence, letting $x \in N$ be any element and taking $X = \{ xg : g
\in S \cup U \}$, we see that $\{ xgf : g, f \in S \cup U \}$ has
support of cardinality $< \kappa$. Also, $\{ x(g + f) : g, f \in S
\cup U \} = \{ xg : g \in S \cup U \} + \{ xf : f \in S \cup U
\}$, as subsets of $N$, and so $\{ x(g + f) : g, f \in S \cup U
\}$ has support of cardinality $< \kappa$. Therefore, by
induction, for all $x \in N$ and $r \in \mathbb{Z}^+$, $B_{S \cup
U} (x, r)$ has support of cardinality $< \kappa$.
\end{proof}

\begin{theorem}\label{fear}
Suppose that $\kappa$ is a regular infinite cardinal $\leq
|\Omega|$, that $M$ can be generated by $< \kappa$ elements as an
$R$-module, and that $S \subseteq E$ is a $\kappa$-fearing
subring. Then $S \not\approx_\kappa E$.
\end{theorem}

\begin{proof}
Let $U \subseteq E$ be a subset such that $|U| < \kappa$.  We will
show that $E \nsubseteq \langle S \cup U \rangle$.

Let $f \in \langle S \cup U \rangle$ be any element.  Then $f$ is
represented by a word of length $r$ in elements of $S \cup U$ for
some $r \in \mathbb{Z}^+$, which implies that $p_{S \cup U}(x, xf)
\leq r$ for every $x \in N$. Hence, if we find an endomorphism $g
\in E$ such that $\{ p_{S \cup U}(x, xg) : x \in N\}$ has no
finite upper bound, then $g \notin \langle S \cup U \rangle$.

In order to construct such a $g$, let us first define two
sequences of elements of $N$. We can pick $x_1 = y_1 \in N$
arbitrarily, and then, assuming that elements with subscripts $i <
j$ have been chosen, we can find a finite subset $\Gamma \subseteq
\Omega$ such that $x_1, \dots, x_{j-1}, y_1, \dots, y_{j-1} \in
M^\Gamma$. Now, let $x_j$ be any nonzero element in $M^{\Omega
\setminus \Gamma}$.  Since $x_j \neq 0$ and $|\Omega \setminus
\Gamma| \geq \kappa$, $M^{\Omega \setminus \Gamma} \cap x_j E$ has
support of cardinality $\geq \kappa$.  (For every $\alpha \in
\Omega \setminus \Gamma$ there is an endomorphism $f \in E$ such
that $y = x_j f$ has the property that $y_\alpha \neq 0$.) Thus,
$(M^{\Omega \setminus \Gamma} \cap x_j E) \setminus B_{S \cup
U}(x_j, j)$ is nonempty, since, by the previous lemma, $B_{S \cup
U}(x_j, j)$ has support of cardinality $< \kappa$; let $y_j$ be
any element thereof.  Now, let $\{\Delta_j : j \in \mathbb{Z}^+\}$
be a collection of disjoint subsets of $\Omega$ such that $x_j,
y_j \in M^{\Delta_j}$ for each $j \in \mathbb{Z}^+$. Let $g_j \in
\E_R(M^{\Delta_j})$ be an endomorphism such that $y_j = x_j g_j$.
Finally, let $g \in E$ be an endomorphism such that the
restriction of $g$ to each $M^{\Delta_j}$ is $g_j$. Such an
endomorphism will have the desired property.
\end{proof}

We note that in the above lemma and theorem the restriction on the
size of $M$ is necessary.  For example, suppose that $M =
\bigoplus_{\aleph_0} L$, for some nonzero left $R$-module $L$, and
also that $\Omega = \aleph_0$. Let $D$ be the diagonal subring of
$\E_R (N)$, consisting of all elements $f \in \E_R (N)$ such that
for each $\alpha \in \Omega$, $M_\alpha f \subseteq M_\alpha$. $D$
is clearly $\aleph_0$-fearing, but $D \approx_{\aleph_0} \E_R
(N)$, since if we take any $f \in \E_R (N)$ that restricts to an
isomorphism $M_\alpha \cong_R N$ for some $\alpha \in \Omega$ and
let $g \in \E_R (N)$ be the inverse of that isomorphism composed
with the inclusion of $M_\alpha$ in $N$, then $\E_R (N) \subseteq
gDf$. In particular, for any nonzero element $x \in N$, $B_{D \cup
\{f, g\}} (x, 3)$ has support of cardinality $\aleph_0$.

In subsequent sections we will focus on the case where $\Omega$ is
countable. However, if $\Omega$ is assumed to be uncountable, then
Lemma~\ref{bowls} can be used to obtain a conclusion stronger than
the one in Theorem~\ref{fear}.

\begin{proposition}
Suppose that $\kappa$ is a regular uncountable cardinal $\leq
|\Omega|$, that $M$ can be generated by $< \kappa$ elements as an
$R$-module, that $S \subseteq E$ is a $\kappa$-fearing subring,
and that $U \subseteq E$ is a subset of cardinality $< \kappa$.
Then $\langle S \cup U \rangle$ is also $\kappa$-fearing.
\end{proposition}

\begin{proof}
Let $x \in N$ be any element. Then, by Lemma~\ref{bowls}, $B_{S
\cup U}(x,r)$ has support of cardinality $< \kappa$ for all $r \in
\mathbb{Z}^+$. As a regular uncountable cardinal, $\kappa$ has
uncountable cofinality, so this implies that $x \langle S \cup U
\rangle (= \bigcup_{r \in \mathbb{Z}^+} B_{S \cup U}(x,r))$ has
support of cardinality $< \kappa$. Also, the support of $R(x
\langle S \cup U \rangle)$ is contained in the support of $x
\langle S \cup U \rangle$, so $(Rx) \langle S \cup U \rangle$ has
support of cardinality $< \kappa$. Now, pick any element $\alpha
\in \Omega$. Letting $x$ range over a generating set of
cardinality $< \kappa$ for $M_\alpha$ as an $R$-module, we
conclude that $M_\alpha \langle S \cup U \rangle$ has support of
cardinality $< \kappa$. Hence, $\langle S \cup U \rangle$ is a
$\kappa$-fearing subring.
\end{proof}

\section{Weakly $\aleph_0$-fearing subrings}

In this section we will keep $R$, $N = \bigoplus_\Omega M$, and
$E$ as before, but will now focus on the case when $\Omega$ is
countable and $M$ is finitely generated. For simplicity, we will
assume that $\Omega = \mathbb{Z}^+$. From now on
$\prec_{\aleph_0,E}$, $\preq_{\aleph_0,E}$, and
$\approx_{\aleph_0,E}$ will be written simply as $\prec$, $\preq$,
and $\approx$, respectively. Also, we will view elements of $E$ as
row-finite matrices over $\E_R(M)$, whenever convenient.

As in the paragraph following Theorem~\ref{fear}, we define $D
\subseteq E$ to be the subring consisting of all elements $f \in
\E_R (N)$ such that for each $\alpha \in \Omega$, $M_\alpha f
\subseteq M_\alpha$. Let $T \subseteq E$ denote the subring of
lower-triangular matrices, consisting of all elements $f \in E$
such that for each $\alpha \in \Omega \ (=\mathbb{Z}^+)$,
$M_\alpha f \subseteq M^\Sigma$, where $\Sigma = \{\gamma \in
\Omega : \gamma \leq \alpha \}$. Also, let $\bar{T} \subseteq E$
denote the subring of upper-triangular matrices, consisting of all
elements $f \in E$ such that for each $\alpha \in \Omega$,
$M_\alpha f \subseteq M^\Gamma$, where $\Gamma = \{\gamma \in
\Omega : \gamma \geq \alpha \}$.

\begin{proposition}\label{D=T}
There exist $g, h \in E$ such that $T \subseteq gDh$. In
particular, $D \approx T$.
\end{proposition}

\begin{proof}
Consider the following two matrices in $E$:
$$A = \left(\begin{array}{c@{}c@{}c@{}c@{}c}
\begin{array}{|c|}\hline
1\\\hline
\end{array} &
\begin{array}{cc}
0 & 0\\
\end{array} &
\begin{array}{ccc}
0 & 0 & 0 \\
\end{array} &
\begin{array}{c}
0 \\
\end{array} &
\begin{array}{c}
\dots \\
\end{array} \\
0 &
\begin{array}{|cc|}\hline
1 & 1\\\hline
\end{array} &
\begin{array}{ccc}
0 & 0 & 0 \\
\end{array} &
\begin{array}{c}
0 \\
\end{array} &
\begin{array}{c}
\dots \\
\end{array} \\
0 &
\begin{array}{cc}
0 & 0\\
\end{array} &
\begin{array}{|ccc|}\hline
1 & 1 & 1 \\\hline
\end{array} &
\begin{array}{c}
0 \\
\end{array} &
\begin{array}{c}
\dots \\
\end{array} \\
\begin{array}{c}
\vdots \\
\end{array} &
\begin{array}{cc}
\vdots & \ \vdots \\
\end{array} &
\begin{array}{ccc}
\vdots & \, \vdots & \ \vdots \\
\end{array} &
\begin{array}{c}
\vdots \\
\end{array} &
\begin{array}{c}
 \ddots \\
\end{array} \\
\end{array}
\right),\ B = \left(\begin{array}{c@{}c}
\begin{array}{c@{}c}
\begin{array}{|c|}\hline
1\\\hline
\end{array} &
\begin{array}{ccc}
0 & 0 & 0 \\
\end{array} \\
\end{array} & \dots \\
\begin{array}{c@{}c}
\begin{array}{|cc|}\hline
1 & 0 \\
0 & 1 \\\hline
\end{array} &
\begin{array}{cc}
0 & 0 \\
0 & 0 \\
\end{array}
\end{array} &
\begin{array}{c}
\dots \\
\dots \\
\end{array}\\
\begin{array}{c@{}c}
\begin{array}{|ccc|}\hline
1 & 0 & 0 \\
0 & 1 & 0 \\
0 & 0 & 1 \\\hline
\end{array} &
\begin{array}{c}
0 \\
0 \\
0 \\
\end{array}
\end{array} &
\begin{array}{c}
\dots \\
\dots \\
\dots \\
\end{array}\\
\begin{array}{cccc}
\vdots & \ \vdots & \ \vdots & \ \vdots \\
\end{array} &
\begin{array}{c}
\ddots \\
\end{array}\\
\end{array}
\right).$$ Let $Y \in T$ be any element. Then we can write
$$Y = \left(\begin{array}{ccccc}
a_{11} & 0 & 0 & 0 & \dots \\
a_{21} & a_{22} & 0 & 0 & \dots \\
a_{31} & a_{32} & a_{33} & 0 & \dots \\
\vdots & \vdots & \vdots & \vdots & \ddots
\end{array}
\right)$$ for some $a_{ij} \in \E_R(M)$. Let $X \in D$ be the
matrix
$$\left(\begin{array}{c@{}c@{}c@{}c}
\begin{array}{|c|}\hline
a_{11}\\\hline
\end{array} &
\begin{array}{ll}
0 & \ \ 0\\
\end{array} &
\begin{array}{lll}
0 & \ \ 0 & \ \ \, 0 \\
\end{array} &
\begin{array}{l}
\dots\\
\end{array} \\
\begin{array}{l}
0\\
0\\
\end{array} &
\begin{array}{|cc|}\hline
a_{21} & 0 \\
0 & a_{22} \\\hline
\end{array} &
\begin{array}{lll}
0 & \ \ 0 & \ \ \, 0 \\
0 & \ \ 0 & \ \ \, 0 \\
\end{array} &
\begin{array}{l}
\dots \\
\dots \\
\end{array}\\
\begin{array}{l}
0 \\
0 \\
0 \\
\end{array} &
\begin{array}{ll}
0 & \ \ 0\\
0 & \ \ 0\\
0 & \ \ 0\\
\end{array} &
\begin{array}{|ccc|}\hline
a_{31} & 0 & 0 \\
0 & a_{32} & 0 \\
0 & 0 & a_{33} \\\hline
\end{array} &
\begin{array}{l}
\dots \\
\dots \\
\dots \\
\end{array} \\
\begin{array}{c}
\vdots \\
\end{array} &
\begin{array}{cc}
\vdots & \ \ \ \vdots \\
\end{array} &
\begin{array}{ccc}
\vdots & \ \ \ \vdots & \ \ \ \vdots \\
\end{array} & \ \ddots \\
\end{array}
\right).$$ Then $AXB = Y$, and so $T \subseteq ADB$. The final
assertion follows from the fact that $D \subseteq T$.
\end{proof}

We note that our definition of $T$ and the above proof make sense
even without assuming that $M$ is finitely generated, since the
elements of $T$ can be viewed as row-finite matrices regardless of
the size of $M$. However, since we are assuming that $M$ is
finitely generated, we have $T \prec E$, by Theorem~\ref{fear}. On
the other hand, we can obtain a result of the opposite sort for
$\bar{T}$.

\begin{corollary}\label{upper-triangular}
$\bar{T} \approx E$.
\end{corollary}

\begin{proof}
By the previous proposition, $T \preq \bar{T}$, since $D \subseteq
\bar{T}$.  Hence, $E = T + \bar{T}$ implies that $\bar{T} \approx
E$.
\end{proof}

Returning to the subring $D$, we can show that all
$\aleph_0$-fearing subrings of $E$ are $\preq D$.  In fact, the
same can be said of a larger class of subrings of $E$.  We will
devote the rest of the section to proving this.

\begin{definition}
We will say that a subring $S \subseteq E$ is \emph{weakly}
$\aleph_0$-\emph{fearing} if $|\{\alpha\in\Omega :
(M_\alpha)S\text{\ has infinite support}\}| < \aleph_0$.
\end{definition}

\begin{lemma}\label{weak fear = fear}
Let $S \subseteq E$ be a weakly $\aleph_0$-fearing subring. Then
there exists an $\aleph_0$-fearing subring $S' \subseteq E$ such
that $S \preq S'$. In particular, $S \not\approx E$.
\end{lemma}

\begin{proof}
Upon enlarging $S$, if necessary, we may assume that $D \subseteq
S$. Let $\Sigma \subseteq \Omega$ be the finite subset consisting
of the elements $\alpha \in \Omega$ such that $(M_\alpha) S$ has
infinite support. Let $S' \subseteq S$ be the subring consisting
of all elements $f \in S$ such that $(M_\alpha)f \subseteq
(M_\alpha)$ for all $\alpha \in \Sigma$, and let $S'' \subseteq S$
be the subring consisting of all elements $f \in S$ such that
$(M_\alpha)f \subseteq (M_\alpha)$ for all $\alpha \notin \Sigma$.
In particular, $S = \pi_{\Omega \setminus \Sigma} S + \pi_{\Sigma}
S \subseteq S' + S''$. Now, $S'$ is $\aleph_0$-fearing; let us
show that $S' \approx S$. To this end, we will demonstrate that
$S'' \subseteq \langle S' \cup U \rangle$ for some finite subset
$U \subseteq E$.

Let $E_{\mathrm{fin}} \subseteq E$ denote the subring consisting
of all elements that have only finitely many off-diagonal entries.
Also, let us write $e_{ij}$ to denote the standard matrix units.
There are countably many such elements, so $\{e_{ij} : i, j \in
\mathbb{Z}^+\} \subseteq U$, for some finite set $U \subseteq E$,
by Theorem~\ref{4-gen}.  Then $S'' \subseteq E_{\mathrm{fin}}
\subseteq \langle D \cup U \rangle$. But, we assumed that $D
\subseteq S'$, so $S'' \subseteq \langle S' \cup U \rangle$, as
desired.

The final assertion follows from Theorem~\ref{fear}.
\end{proof}

\begin{lemma}\label{fear <T}
Let $S \subseteq E$ be an $\aleph_0$-fearing subring.  Then $S
\subseteq gT$ for some $g \in E$.
\end{lemma}

\begin{proof}
For each $k \in \Omega \, (=\mathbb{Z}^+)$, let $l_k$ be the
largest element in the union of the supports of $(M_j)S$  for all
$j \leq k$. Let $f \in E$ be the endomorphism that takes $M_{l_k}$
identically to $M_k$ for each positive integer $k$, and takes
$M_\alpha$ to zero if $\alpha \neq l_k$ for all $k \in
\mathbb{Z}^+$.  Also, let $g \in E$ be the endomorphism that takes
$M_k$ identically to $M_{l_k}$ for each positive integer $k$. Now,
let $h \in S$ be any element. Then $fh \in T$, and $g(fh) = (gf)h
= 1 \cdot h = h$. Hence, $S \subseteq gT$.
\end{proof}

As in Proposition~\ref{D=T}, the above proof works even without
assuming that $M$ is finitely generated.

The following theorem summarizes the results of this section.

\begin{theorem}\label{weak fear < D}
Let $S \subseteq E$ be a weakly $\aleph_0$-fearing subring. Then
$S \preq D$.
\end{theorem}

\begin{proof}
This follows from Proposition~\ref{D=T}, Lemma~\ref{weak fear =
fear}, and Lemma~\ref{fear <T}.
\end{proof}

\section{Subrings arising from preorders}

Keeping the notation from the previous section, we will now turn
our attention to subrings $S \subseteq E$ such that $D \preq S$.
For this, we will need a new concept.

\begin{definition}\label{full}
Let $\rho$ be a preordering of $\Omega$. Define $E(\rho) \subseteq
E$ to be the subset consisting of those elements $f \in E$ such
that for all $\alpha, \beta \in \Omega$, $\pi_\alpha f \pi_\beta
\neq 0$ implies $(\alpha, \beta) \in \rho$.
\end{definition}

It is clear that the subsets $E(\rho)$ are subrings. For example,
$D$, $T$, $\bar{T}$, and $E$ are of this form. Indeed, recalling
that $\Omega = \mathbb{Z}^+$, and setting $\rho_1 = \{(\alpha,
\alpha) \in \Omega \times \Omega\}$, $\rho_2 = \{(\alpha, \beta)
\in \Omega \times \Omega : \alpha \geq \beta\}$, and $\rho_3 =
\{(\alpha, \beta) \in \Omega \times \Omega : \alpha \leq \beta\}$,
we have $D = E(\rho_1)$, $T = E(\rho_2)$, and $\bar T =
E(\rho_3)$. We also note that every subring of $E$ of the form
$E(\rho)$ contains $D$ and is closed in the function topology
(i.e., the topology inherited from the set $N^N$ of all functions
$\bigoplus_\Omega M \to \bigoplus_\Omega M$, where a subbasis of
open sets is given by the sets $\{f \in N^N : mf= n\}$, for all
$m, n \in \bigoplus_\Omega M$). In fact, if $\E_R(M)$ is a simple
ring, then this characterizes such subrings of $E$.

\begin{proposition}\label{simple full}
Suppose that $\E_R(M)$ is a simple ring, and let $S \subseteq E$
be a subring. Then $S = E(\rho)$ for some preordering $\rho$ of
$\Omega$ if and only if $S$ is closed in the function topology and
$D \subseteq S$.
\end{proposition}

\begin{proof}
Suppose that $S$ is closed in the function topology and $D
\subseteq S$. Let $\rho = \{(\alpha, \beta) : \pi_\alpha S
\pi_\beta \neq 0 \} \subseteq \Omega \times \Omega$. Since $S$
contains the identity element, $\rho$ is reflexive.

Next, we note that for all $\alpha, \beta \in \Omega$ there is an
obvious bijection between $\pi_\alpha E \pi_\beta$ and $\E_R(M)$,
under which $\pi_\alpha S \pi_\beta$ corresponds to a 2-sided
ideal (since $D \subseteq S$). Hence, $\E_R(M)$ being simple
implies that either $\pi_\alpha S \pi_\beta = 0$ or $\pi_\alpha S
\pi_\beta = \pi_\alpha E \pi_\beta$. In particular, if $\pi_\alpha
S \pi_\beta \neq 0$ and $\pi_\beta S \pi_\gamma \neq 0$ for some
$\alpha, \beta, \gamma \in \Omega$, then $\pi_\alpha S \pi_\beta
\pi_\beta S \pi_\gamma = \pi_\alpha E \pi_\beta \pi_\beta E
\pi_\gamma = \pi_\alpha E \pi_\gamma$. Since $\pi_\beta \in D
\subseteq S$, we have $0 \neq \pi_\alpha S \pi_\beta \pi_\beta S
\pi_\gamma \subseteq \pi_\alpha S \pi_\gamma$. Hence, $\rho$ is
transitive and therefore a preorder.

Let $f \in E$ be an element with the property that $\pi_\alpha f
\pi_\beta \neq 0$ implies $(\alpha, \beta) \in \rho$. Then
$\pi_\alpha f \pi_\beta \in \pi_\alpha E \pi_\beta \subseteq S$
for all $\alpha, \beta \in \Omega$, by the previous paragraph.
Now, $f$ is in the closure of the set of sums of elements of the
form $\pi_\alpha f \pi_\beta$. Hence $f \in S$, by the hypothesis
that $S$ is closed in the function topology. This shows that $S =
E(\rho)$.

The converse is clear.
\end{proof}

In the previous section we showed that if $M$ is finitely
generated and $S \subseteq E$ is a weakly $\aleph_0$-fearing
subring, then $S \not\approx E$.  It turns out that for a subring
$S \subseteq E$ of the form $E(\rho)$ the property of being weakly
$\aleph_0$-fearing is not only sufficient but also necessary for
$S \not\approx E$. We require a lemma before proceeding to the
proof of this statement.

\begin{lemma}\label{sets}
Let $\Phi$ and $\Gamma$ be sets such that $|\Phi| = |\Gamma| =
\aleph_0$, and let $\{\Lambda_\varphi : \varphi \in \Phi \}$ be a
collection of infinite subsets of $\Gamma$. Then there is a subset
$\{\varphi_j : j \in \mathbb{Z}^+ \} \subseteq \Phi$ of distinct
elements and a collection of infinite sets $\{\bar{\Lambda}_j : j
\in \mathbb{Z}^+ \}$, where for each $j \in \mathbb{Z}^+$,
$\bar{\Lambda}_j \subseteq \Lambda_{\varphi_j}$, such that one of
the following holds:
\begin{enumerate}
\item[$(1)$] $\bar{\Lambda}_j \subseteq
\bar{\Lambda}_{j'}$ for $j \geq j'$,
\item[$(2)$] $\bar{\Lambda}_j \cap
\bar{\Lambda}_{j'} = \emptyset$ for $j \neq j'$.
\end{enumerate}
\end{lemma}

\begin{proof}
Suppose that $\{\varphi_j : j \in \mathbb{Z}^+ \}$ is a sequence
of distinct elements of $\Phi$. Let us construct inductively an
infinite collection of sets $\{\bar{\Lambda}_j : j \in
\mathbb{Z}^+ \}$, where for each $j \in \mathbb{Z}^+$,
$\bar{\Lambda}_j \subseteq \Lambda_{\varphi_j}$, and
$\bar{\Lambda}_{j+1} \subseteq \bar{\Lambda}_j$. Let
$\bar{\Lambda}_1 = \Lambda_{\varphi_1}$, and assuming that
$\bar{\Lambda}_n$ has been constructed, let $\bar{\Lambda}_{n+1} =
\bar{\Lambda}_n \cap \Lambda_{\varphi_{n+1}}$.

If there is a sequence of distinct elements of $\Phi$,
$\{\varphi_j : j \in \mathbb{Z}^+ \}$, such that the collection
constructed above consists of infinite sets, then (1) is
satisfied. So suppose that no such sequence exists. Let
$\varphi_1, \varphi_2, \dots, \varphi_n \in \Phi$ be any sequence
of elements such that $n \geq 1$, $\Delta_1 := \Lambda_{\varphi_1}
\cap \Lambda_{\varphi_2} \cap \dots \cap \Lambda_{\varphi_n}$ is
infinite, and for all $\varphi \in \Phi \setminus \{\varphi_1,
\varphi_2, \dots, \varphi_n\}$, $\Delta_1 \cap \Lambda_\varphi$ is
finite. Set $\Phi_1 = \Phi \setminus \{\varphi_1, \varphi_2,
\dots, \varphi_n\}$. Repeating the above process, let $\varphi_1',
\varphi_2', \dots, \varphi_n' \in \Phi_1$ be any sequence of
elements such that $n \geq 1$, $\Delta_2 := \Lambda_{\varphi_1'}
\cap \Lambda_{\varphi_2'} \cap \dots \cap \Lambda_{\varphi_n'}$ is
infinite, and for all $\varphi \in \Phi_1 \setminus \{\varphi_1',
\varphi_2', \dots, \varphi_n'\}$, $\Delta_2 \cap \Lambda_\varphi$
is finite. Set $\Phi_2 = \Phi_1 \setminus \{\varphi_1',
\varphi_2', \dots, \varphi_n'\}$, etc. Continuing in this fashion,
we obtain a subset $\{\varphi_j : j \in \mathbb{Z}^+ \} \subseteq
\Phi$ of distinct elements and a collection of infinite sets
$\{\Delta_j : j \in \mathbb{Z}^+ \}$, where for each $j \in
\mathbb{Z}^+$, $\Delta_j \subseteq \Lambda_{\varphi_j}$, and
$\Delta_j \cap \Delta_{j'}$ is finite for $j \neq j'$.

Let $\bar{\Lambda}_1 = \Delta_1$, and for each $j > 1$, let
$\bar{\Lambda}_j = \Delta_j \setminus \bigcup_{i=1}^{j-1}
(\Delta_i \cap \Delta_j)$.  Then the elements of
$\{\bar{\Lambda}_j : j \in \mathbb{Z}^+ \}$ satisfy (2).
\end{proof}

\begin{proposition}\label{full+weak = E}
Suppose that $M$ is finitely generated, $\rho$ is a preordering of
$\Omega$, and $S = E(\rho)$. If $S$ is not weakly
$\aleph_0$-fearing, then $S \approx E$.
\end{proposition}

\begin{proof}
For each $\alpha \in \Omega \ (=\mathbb{Z}^+)$, denote the support
of $(M_\alpha)S$ by $\mathrm{supp}(M_\alpha S)$, and consider the
set $\{\mathrm{supp}(M_\alpha S) : |\mathrm{supp}(M_\alpha S)| =
\aleph_0\}$. This is an infinite collection of infinite subsets of
$\Omega$, since $S$ is not weakly $\aleph_0$-fearing. By the
previous lemma, there is a set $\Sigma = \{\alpha_j : j \in
\mathbb{Z}^+\}$ of distinct elements of $\Omega$ and a collection
of infinite sets $\{\bar{\Lambda}_j : j \in \mathbb{Z}^+ \}$,
where for each $j \in \mathbb{Z}^+$, $\bar{\Lambda}_j \subseteq
\mathrm{supp}(M_{\alpha_j} S)$, and either
\begin{enumerate}
\item[(1)] $\bar{\Lambda}_j \subseteq \bar{\Lambda}_{j'}$ for
$j \geq j'$, or
\item[(2)] $\bar{\Lambda}_j \cap \bar{\Lambda}_{j'} = \emptyset$ for $j
\neq j'$.
\end{enumerate}
We will now treat the two cases individually.

Suppose that (1) holds.  We begin by constructing a sequence
$\{\beta_j : j \in \mathbb{Z}^+\}$ of elements of
$\bar{\Lambda}_1$. Pick $\beta_1 \in \bar{\Lambda}_1$ arbitrarily.
Let $\beta_2 \in \bar{\Lambda}_2$ be such that $\beta_2 \neq
\beta_1$. Then let $\beta_3 \in \bar{\Lambda}_3$ be such that
$\beta_3 \neq \beta_2$ and $\beta_3 \neq \beta_1$, and so on. Now,
let $S' \subseteq S$ be the subset consisting of all endomorphisms
$f$ such that for each $j \in \mathbb{Z}^+$, $(M_{\alpha_j})f
\subseteq M^{\Gamma_j}$, where $\Gamma_j = \{\beta_i : i \geq j\}
\subseteq \bar{\Lambda}_j$. Let $g \in E$ be the endomorphism that
maps $N$ to $M^\Sigma$ by sending $M_j$ identically to
$M_{\alpha_j}$ for each $j \in \mathbb{Z}^+ = \Omega$.  Also, let
$h \in E$ be an endomorphism that takes $M^{\Gamma_1}$ to $N$ by
sending $M_{\beta_j}$ identically to $M_j$ for each $j \in
\mathbb{Z}^+ = \Omega$. Then $\bar{T} = gS'h$, since $S =
E(\rho)$. Hence, by Corollary~\ref{upper-triangular}, $S \approx
E$.

Suppose that (2) holds.  Then $\sum_{j \in \mathbb{Z}^+}
M^{\bar{\Lambda}_j}$ is direct, and so there is an endomorphism $h
\in E$ that simultaneously maps each $M^{\bar{\Lambda}_j}$
isomorphically to $N$. Let $S' \subseteq S$ be the subset
consisting of all endomorphisms $f$ such that for each $j \in
\mathbb{Z}^+$, $(M_{\alpha_j})f \subseteq M^{\bar{\Lambda}_j}$.
Then $\pi_\Sigma \mathrm{Hom}_R(M^\Sigma, N) \subseteq S'h$, since
$S = E(\rho)$. Now, let $g \in E$ be an endomorphism that maps $N$
isomorphically to $M^\Sigma \subseteq N$.  Then $E \subseteq
gS'h$, and hence $S \approx E$.
\end{proof}

\begin{corollary}\label{E-class}
Suppose that $M$ is finitely generated, and let $\rho$ be a
preordering of $\Omega$. Then $E(\rho) \approx E$ if and only if
$E(\rho)$ is not weakly $\aleph_0$-fearing.
\end{corollary}

\begin{proof}
The forward implication was proved in Lemma~\ref{weak fear =
fear}, while the backward implication was proved in
Proposition~\ref{full+weak = E}.
\end{proof}

Putting together Corollary~\ref{E-class}, Theorem~\ref{weak fear <
D}, and the remarks after Definition~\ref{full} we obtain the
following result.

\begin{theorem}\label{classification of preorder
subrings}
Suppose that $M$ is finitely generated, and let $\rho$
be a preordering of $\Omega$. Then exactly one of the following
holds:
\begin{enumerate}
\item[$(1)$] $E(\rho) \approx D,$
\item[$(2)$] $E(\rho) \approx E.$
\end{enumerate}
\end{theorem}

If we assume that $\E_R(M)$ is a simple ring, then this theorem
can be stated without reference to preorders.

\begin{corollary}\label{classification for simples}
Suppose that $M$ is finitely generated and $\E_R(M)$ is a simple
ring. If $S \subseteq E$ is a subring that is closed in the
function topology and $D \subseteq S$, then exactly one of the
following holds:
\begin{enumerate}
\item[$(1)$] $S \approx D,$
\item[$(2)$] $S \approx E.$
\end{enumerate}
\end{corollary}

\begin{proof}
This follows from Theorem~\ref{classification of preorder
subrings} and Proposition~\ref{simple full}.
\end{proof}

It would be desirable to relax the condition $D \subseteq S$ in
the above statement, say, by instead considering closed subrings
$S$ satisfying $I \subseteq S$, where $I$ is the diagonally
embedded copy of $\E_R(M)$ in $E$. However, doing so makes the
situation much more messy, though we can show that in general such
subrings fall into at least four equivalence classes. In the next
result, let $C \subseteq E$ denote the subring $\langle I \cup
\mathrm{Hom}_R(N, M_1) \iota_1 \rangle = I + \mathrm{Hom}_R(N,
M_1) \iota_1$, where $\iota_1$ is the inclusion of $M_1$ in $N$.

\begin{proposition}\label{separate 4 classes}
Suppose that $M$ is simple and $\E_R (M)$ has countable dimension
as a vector space over its center. Then $I \prec C \prec D \prec
E$.
\end{proposition}

\begin{proof}
By Theorem~\ref{fear}, $D \prec E$, and, by Theorem~\ref{weak fear
< D}, $C \preq D$. Also, $I \subseteq C$. Thus, it suffices to
show that $I \not\approx C \not\approx D$.

Let $Z$ denote the center of $\E_R (M)$.  Also, let $U \subseteq
E$ be a finite subset. Then $\langle I \cup U \rangle$ has
countable dimension as a vector space over $Z$, and hence $C
\not\subseteq \langle I \cup U \rangle$, since $C$ has uncountable
dimension. Therefore $I \not\approx C$.

It remains to show that $C \not\approx D$. By
Proposition~\ref{D=T}, it suffices to prove that given a finite
set $U \subseteq E$, we have $T \not\subseteq \langle C \cup U
\rangle$. Now, for any subset $\Sigma \subseteq \Omega$, let
$\iota_\Sigma$ denote the inclusion of $M^\Sigma$ in $N$, and set
$F = \{h \in E : \exists \Sigma \subseteq \Omega \
\mathrm{finite}, \ \mathrm{such} \ \mathrm{that} \ h \in
\mathrm{Hom}_R(N, M^\Sigma) \iota_\Sigma \}$ (i.e., the set of
matrices with zeros in all but finitely many columns). We will
first show that any $f \in \langle C \cup U \rangle$ can be
written as $f = g + h$, where $g \in \langle I \cup U \rangle$ and
$h \in F$.

Setting $H = \mathrm{Hom}_R(N, M_1) \iota_1$, we have $C = I + H$.
Hence, every element $f \in \langle C \cup U \rangle$ can be
expressed as $f = g + \bar{h}$, where $g \in \langle I \cup U
\rangle$ and $\bar{h}$ is a sum of products of elements of $I$,
$H$, and $U$, such that each product contains an element of $H$.
Now, $EH \subseteq H$, so $\bar{h}$ is a sum of products of the
form $h'g'$, where $h' \in H$ and $g' \in \langle I \cup U
\rangle$. But, such elements $h'g'$ belong to $F$, and hence $f$
can be written as $f = g + h$, where $g \in \langle I \cup U
\rangle$ and $h \in F$ (since $F$ is closed under addition).

Now, for each $f \in T$ pick some $h_f \in F$.  Then $\{f - h_f :
f \in T\}$ generates a $Z$-vector space of uncountable dimension
and is therefore not contained in any subring of $E$ that has
countable dimension, regardless of how the elements $h_f$ are
picked. In particular, $\{f - h_f : f \in T\} \not\subseteq
\langle I \cup U \rangle$. Hence $T \not\subseteq \langle C \cup U
\rangle$, by our description of the elements of $\langle C \cup U
\rangle$ above, completing the proof.
\end{proof}

\noindent
Department of Mathematics \\
University of California \\
Berkeley, CA 94720 \\
USA \\

\noindent Email: {\tt zak@math.berkeley.edu}
\end{document}